\documentclass{crmp-l}

\newtheorem{theorem}{Theorem}[section]
\newtheorem{lemma}[theorem]{Lemma}
\newtheorem{proposition}[theorem]{Proposition}
\newtheorem{corollary}[theorem]{Corollary}
\newtheorem{conjecture}[theorem]{Conjecture}

\theoremstyle{definition}
\newtheorem{definition}[theorem]{Definition}
\newtheorem{example}[theorem]{Example}

\theoremstyle{remark}
\newtheorem{remark}[theorem]{Remark}

\font\tenmsa=msam10
\font\sevenmsa=msam7
\font\fivemsa=msam5
\newfam\msafam
 \textfont\msafam=\tenmsa  \scriptfont\msafam=\sevenmsa
 \scriptscriptfont\msafam=\fivemsa
\def\hexnumber@#1{\ifcase#1 0\or1\or2\or3\or4\or5\or6\or7\or8\or9\or
A\or B\or C\or D\or E\or F\fi }
\edef\msa@{\hexnumber@\msafam}
\mathchardef\circledS="0\msa@73

\numberwithin{equation}{section}

\begin{document}

\title[Picard-Fuchs Uniformization]{Picard-Fuchs Uniformization: Modularity of the Mirror Map and Mirror-Moonshine}

\author{Charles F. Doran}
\address{Department of Mathematics, Harvard University, Cambridge, Massachusetts 02138}
\email{doran@math.harvard.edu}

\date{November 20, 1998}

\begin{abstract}

Motivated by a conjecture of Lian and Yau concerning the mirror map in string theory \cite{LY1,ICMP}, we determine when the mirror map $q$-series of certain elliptic curve and K3 surface families are Hauptmoduln (genus zero modular functions).  Our geometric criterion for modularity characterizes orbifold uniformization properties of their Picard-Fuchs equations, effectively demystifying the {\em mirror-moonshine phenomenon}.
A longer, more comprehensive treatment of these results can be found in \cite{Dor2}.  For a detailed look at several explicit examples of this phenomenon, see the article by Verrill and Yui in this volume \cite{VY}.  
\end{abstract}

\maketitle

\tableofcontents

\section{Introduction}

Although the {\em mirror map} was introduced into string theoretic physics only in the last decade, whence it has lead to mathematical predictions of curve counts on certain Calabi-Yau threefolds, it is in fact reminiscent of standard mathematical constructions of the 19th century dating back at least to Legendre's work on elliptic functions.  In a similarly classical vein such maps have long been of mathematical interest through their connection with modular functions, e.g., the Hauptmoduln of Fricke-Klein.  Stimulated by the Mirror-Moonshine Conjecture of Lian-Yau (\cite{LY1,ICMP}, Conjecture \ref{LYConj} below), in this paper we characterize families of elliptic curves and K3 surfaces for which the mirror map is {\em modular}, i.e., is a Hauptmodul.

The Lian-Yau conjecture was originally stated just for a small class of K3 surface pencils constructed by an ``orbifold construction'', but it was rapidly realized that their {\em mirror-moonshine} observations apply to a larger class of ``torically defined'' K3 surface families \cite{LY2}.  In fact, as we show here, the modularity phenomenon they observed can occur quite generally in families of elliptic curves with section and families of lattice polarized K3 surfaces.  In particular, the phenomenon is not ``toric'' in origin, and can be characterized by other means.  The strategy we follow to do so makes use of the existence of a one-dimensional (coarse) moduli space and a global Torelli theorem for each type of family.  In the elliptic curve case the moduli space is just the $J$-line and the Torelli theorem is classical.  For families of K3 surfaces with a lattice polarization such moduli and Torelli results were recently formalized by Dolgachev for his own investigations of mirror-symmetric families of K3 surfaces motivated by Arnold's ``strange duality'' \cite{Dol}.

We begin in \S \ref{seccy} by introducing Calabi-Yau manifolds, Picard-Fuchs equations, and the mirror map.  The periods of the holomorphic top forms for a family of Calabi-Yau manifolds define multivalued solutions to the associated {\em Picard-Fuchs differential equation}, a homogeneous Fuchsian ordinary differential equation on the base of the family.  The {\em mirror map} is then a ``local inverse'' to the projectivized period mapping at a regular singular point with maximal unipotent monodromy.  In fact, since the mirror map depends only on the {\em projective} periods, it is insensitive to simultaneous rescalings of these periods.  Thus it is really the {\em projective normal form} of the Picard-Fuchs differential equation which determines the mirror map.  

Second order homogeneous Fuchsian ordinary differential equations in projective normal form have a particular connection with automorphic functions, dating back at least to Poincar\'{e} and Fricke-Klein.  An automorphic function $f$ for a genus zero Fuchsian group of the first kind $\Gamma$ can always be realized as the single-valued inverse of the ratio of two multi-valued functions $u_1(x)$ and $u_2(x)$ on the curve $\Gamma \, \backslash \, {\mathbb{H}}^* \simeq {\mathbb{P}}^1$.  The functions $u_i$ satisfy a second order Fuchsian ordinary differential equation in projective normal form which is naturally associated to $f$ (\cite{Leh} and \S \ref{ssorbuni} below). This equation is sometimes called the ``Schwarzian equation'' or {\em (orbifold) uniformizing differential equation} attached to the automorphic function.  We conclude that, at least for families of Calabi-Yau manifolds with order two Picard-Fuchs equations:
\begin{quotation}
The mirror map is an automorphic function if and only if the projective normal form of the Picard-Fuchs differential equation is an orbifold uniformizing differential equation.  
\end{quotation}

In fact, families of elliptic curves with section do have second order Picard-Fuchs equations (\S \ref{PFWei}).  For the families of $M_n$-polarized K3 surfaces, Picard-Fuchs equations are of order three, but they take a particular form as the {\em symmetric square} of an associated second order equation, a process which preserves the period ratio and hence the mirror map (\S \ref{PFMn}).  Thus, to solve the mirror-moonshine modularity problem we need only characterize, in terms of a ``geometric'' invariant of the family of elliptic curves with section (respectively $M_n$-polarized K3 surfaces), when the projective normal form of the Picard-Fuchs equation (respectively its {\em square root}) is a uniformizing differential equation attached to an automorphic function.  We call this {\em Picard-Fuchs uniformization}.

In the case of elliptic curves, the natural geometric invariant to consider is Kodaira's {\em functional invariant} $\mathcal{J}(z)$, defined as the composition of the multi-valued period morphism $\tau(z) = \omega_1(z) / \omega_0(z)$ to the upper half plane and the elliptic modular function $J(\tau)$ to the $J$-line.  In other words, the functional invariant associates to each point of the base of the family the $J$-invariant of its elliptic curve fiber.  Kodaira proves that $\mathcal{J}(z)$ is a rational function.  We consider the ``pullback by $\mathcal{J}$'' of the uniformizing differential equation $\Lambda$ for $PSL(2,\mathbb{Z})$.  By the rationality of $\mathcal{J}(z)$, the resulting differential equation $\mathcal{J}^*(\Lambda)$ agrees with the Picard-Fuchs equation up to projective equivalence (Proposition \ref{PropPFLJ}).  Thus we need only determine the conditions on the rational function $\mathcal{J}(z)$ for the projective normal form of $\mathcal{J}^*(\Lambda)$ to be a uniformizing differential equation.  We accomplish this by applying our general criterion, Theorem \ref{thpullcharexp}, which characterizes rational functions that pullback one projective equivalence class of uniformizing differential equations to another.

Next, in \S \ref{Mnpol}, we extend our criterion to the case of $M_n$-polarized families of K3 surfaces.  Dolgachev has shown that the (coarse) moduli space for $M_n$-polarized K3 surfaces is uniformized by the group $\Gamma_0(n)^*$ obtained by adding the Fricke involution to the congruence subgroup $\Gamma_0(n)$ \cite{Dol}.  By analogy with the elliptic curve case, we introduce the {\em generalized functional invariant} $\mathcal{H}_n(z)$ as the composition of the multi-valued truncated period morphism $\tau(z) = \Omega_1(z) / \Omega_0(z)$ to the upper half plane and the normalized Hauptmodul $H_n(\tau)$ for $\Gamma_0(n)^*$, and establish its rationality (\S \ref{genfct}).  The criterion from \S \ref{regsingp} then allows us to determine the conditions on $\mathcal{H}_n(z)$ to pullback the uniformizing differential equation for $\Gamma_0(n)^*$ to another uniformizing differential equation up to projective equivalence.  As a consequence we obtain the desired characterization of the Lian-Yau mirror-moonshine phenomenon in Theorem \ref{LYchar}.

Finally, in \S \ref{general}, a number of natural generalizations of this work are discussed.  Our Picard-Fuchs uniformization results for $M_n$-polarized K3 surface families extend {\em mutatis mutandis} to the general rank 19 lattice polarized case.  We note that the genus-zero restriction on both the uniformizing groups and the base of the families of elliptic curves and K3 surfaces involved can be removed.  Moreover, there is an extension of our modularity criteria to certain multi-parameter families of abelian varieties and K3 surfaces of various Picard ranks.  In fact, there are one-parameter families of Picard rank 18 K3 surfaces whose mirror maps are ``bi-modular''.  These generalizations are treated in \cite{Dor2}.

Throughout this paper we assume some basic familiarity with the theory of Fuchsian ordinary differential equations: regular singular points and their classification into types (i.e., generic, logarithmic, apparent), characteristic exponents, Frobenius' theorem, etc.  A favorite reference is Yoshida's book \cite{Yos} which also includes the theory of (orbifold) uniformization.  Other nice presentations of the classical Fuchsian theory, with quite different applications, can be found in the survey paper by Varadarajan \cite{Var} and the book by Anosov and Bolibruch \cite{AB}.

\begin{remark}
This article is intended to complement the explicit constructions described by Helena Verrill and Noriko Yui in this proceedings \cite{VY}.  By providing a simple, geometric criterion for the Lian-Yau mirror-moonshine phenomenon which applies equally well to both elliptic curve and K3 surface families, we hope at once to ``demystify'' such examples of this phenomenon, and to highlight the importance of explicit constructions in answering questions which remain.  Taking advantage of the opportunity such a pairing affords, we refer the reader to their paper for several explicit examples and extra discussion of related subject matter.
\end{remark}

\section{Calabi-Yau manifolds and the mirror map}

\label{seccy}

In this section we quickly review the basics of Calabi-Yau manifolds, Picard-Fuchs equations, and the mirror map.  The definitions and concepts are stated without proof, and a single motivational example is provided: that of the universal family of elliptic curves over the $J$-line.  Perhaps the chief novelty here is the emphasis on the role played by the ``projective normal form'' Fuchsian ordinary differential equation.  Although our discussion applies quite generally to families of Calabi-Yau varieties, we will only be applying it in this paper to the elliptic curve and K3 surface cases.  

At the time of writing, there are only a few general surveys on Calabi-Yau manifolds and the mirror map written for mathematicians, e.g., \cite{EMM, MSI}, \cite{BP}, \cite{Voi}.  By the time this volume goes to print it is hoped that the compendium \cite{CK} will have appeared to fill this gap in the literature.

\subsection{Calabi-Yau manifolds}

\mbox{}

\begin{definition}
A {\em Calabi-Yau manifold} $M$ is a compact complex manifold with trivial canonical bundle.
\end{definition}

\begin{example}
A one dimensional Calabi-Yau manifold is an elliptic curve.  A simply connected Calabi-Yau manifold of dimension two is a K3 surface.
\end{example}

Up to complex scaling, a Calabi-Yau $n$-fold has a unique holomorphic $n$-form 
$$\omega_{n,0} \in H^n(M, {\mathbb{C}}) \ .$$

\begin{definition}
The {\em periods} of this holomorphic $n$-form are the numbers
$$p_\gamma(M) := \int_\gamma \omega_{n,0}$$
as $\gamma$ runs through the $n$-cycles of $M$.  Fix an ordered basis of $n$-cycles $\{ \gamma_1, \ldots, \gamma_{k} \}$, for $k$ the $n$th Betti number of $M$.  Then the periods of $M$ define the vector
$$[p_{\gamma_i}(M)]_{i=1}^k \in {\mathbb{P}}^{k-1}$$
the (projective) {\em period point} of $M$.
\end{definition}

\begin{example}
Consider the elliptic curve $y^2 = 4 x^3 - g_2 x - g_3$ in Weierstrass form.  Here $n = 1$ and $k = 2$.  The holomorphic one form is the usual differential of the first kind $dx / y$. 
\end{example}

\subsection{Period mapping and Picard-Fuchs equations}

\mbox{}

\begin{definition}
Given a family $\pi : X \rightarrow S$ of Calabi-Yau manifolds $X_s$, by continuous extension of the basis of $n$-cycles the periods of the Calabi-Yau fibers define a vector of $k$ multi-valued functions on the base manifold $S$
$$[p_{\gamma_i(s)}(X_s)]_{i=1}^k : S \rightarrow {\mathbb{P}}^{k-1}$$
the {\em period mapping} associated with the family.
\end{definition}

\begin{example}
Consider the family $\mathcal{E}$ of elliptic curves over ${\mathbb{P}}^1$ defined by the equation
$$ \mathcal{E} : y^2 = 4 x^3 - \frac{27 s}{s - 1} x - \frac{27 s}{s - 1} \ .$$
The periods of the form $dx/y$ may be given in terms of the hypergeometric function $_2F_1$ (see \cite{Sti}, pp.~232--233, for explicit expressions).
\label{egellJ}
\end{example}

Suppose now that $S = {\mathbb{P}}^1$.

The functions $p_{\gamma(s)}(X_s)$ satisfy an ordinary differential equation with regular singular points (i.e., a {\em Fuchsian} ODE) \cite{BP}
\begin{eqnarray}
\frac{d^k f}{ds^k} + P_1(s) \frac{d^{k-1} f}{ds^{k-1}} + \ldots + P_k(s) f &=& 0 \ , \ P_i(s) \in {\mathbb{C}}(s) \ .
\label{eqnPF}
\end{eqnarray}
\begin{definition}
The differential equation (\ref{eqnPF}) is called the {\em Picard-Fuchs differential equation} of the family $\pi : X \rightarrow {\mathbb{P}}^1$.  
\end{definition}

\begin{example}
For the family $\mathcal{E}$ in Example \ref{egellJ}, the Picard-Fuchs equation is
$$\frac{d^2 f}{ds^2} + \frac{1}{s} \frac{df}{ds} + \frac{(31/144) s - 1/36}{s^2 (s-1)^2} f = 0 \ .$$
There is a basis of solutions with local monodromies $G_0$, $G_1$, $G_\infty$ about the regular singular points $\{0, 1, \infty\}$ respectively, where
$$G_0 = \left( \begin{array}{cc} 1 & 1 \\ -1 & 0 \end{array} \right) \ , \ G_1 = \left( \begin{array}{cc} 0 & -1 \\ 1 & 0 \end{array} \right) \ , \ G_\infty = \left( \begin{array}{cc} 1 & 1 \\ 0 & 1 \end{array} \right) \ .$$
\label{egPFJ}
\end{example}

\subsection{Points of maximal unipotent monodromy and the mirror map}

\mbox{}

\begin{definition}
A regular singular point of a Fuchsian ordinary differential equation is called a {\em point of maximal unipotent monodromy} if the local monodromy matrix $G$ is such that $G - I_k$ is nilpotent with exact order $k$.
\end{definition}
In a neighborhood of a point of maximal unipotent monodromy, Frobenius' method tells us that there is a basis of solutions such that the first is holomorphic at the point, the second has logarithmic behavior, the next behaves like $\log^2$, $\ldots$, up to $\log^{k-1}$.

\begin{example}
 The point of maximal unipotent monodromy in Example \ref{egPFJ}  is  $\infty \in {\mathbb{P}}^1$.
\end{example}

\begin{remark}
Not every family of Calabi-Yau manifolds will have a Picard-Fuchs differential equation with a point of maximal unipotent monodromy.  For a class of {\em torically constructed} hypersurfaces in toric Fano varieties, the existence of points of maximal unipotent monodromy was established by Hosono, Lian, and Yau \cite{HLY}.
\end{remark}

Consider a family of Calabi-Yau manifolds $\pi : X \rightarrow {\mathbb{P}}^1$, whose Picard-Fuchs equation has a point of maximal unipotent monodromy.  In a neighborhood of such a point consider the truncated period vector consisting only of the holomorphic solution and the logarithmic solution
$$[p_{hol}(s), p_{log}(s)] : {\mathbb{P}}^1 \rightarrow {\mathbb{P}}^1 \ .$$
If the image lies in the upper half plane ${\mathbb{H}} \subset {\mathbb{P}}^1$, then, possibly after composition with projective linear transformations so that the singular point lies at $0 \in {\mathbb{P}}^1$ and maps to $\imath \infty \in {\mathbb{H}}^* \subset {\mathbb{P}}^1$, we can consider the $q$-series for the local inverse mapping
$$z(q(\tau)) : {\mathbb{H}} \rightarrow {\mathbb{P}}^1 \ , \  q(\tau) = e^{2 \pi \imath \tau} \ .$$
\begin{definition}
This $q$-series $z(q)$ is called the {\em mirror map} of the family $\pi : X \rightarrow {\mathbb{P}}^1$ about the point of maximal unipotent monodromy.
\end{definition}

\begin{example}
 For the family $\mathcal{E}$ of Example \ref{egellJ} the mirror map is quite familiar.  Since the maximal unipotent monodromy point is at $\infty$, we change variables first to $z = 1/s$.  The single-valued local inverse to the period mapping  is then the reciprocal of the $q$-series for the elliptic modular function $J(q)$
$$J(q) = \frac{1}{q} + 744 + 196884 q + 21493760 q^2 + O(q^3) \ , $$
$$z(q) = \frac{1}{J(q)} = q - 744 q^2 + 356652 q^3 - 140361152 q^4 + O(q^5) \ .$$
\end{example}

\begin{remark}
The period mapping is defined as a map to projective space.  If one is interested in the mirror map it is often preferable to consider the Picard-Fuchs differential equation only up to ``projective equivalence''.  
\begin{definition}
The {\em projective normal form} of a Fuchsian ordinary differential equation (e.g., that in equation (\ref{eqnPF}) above) is the unique Fuchsian ordinary differential equation without a $(k-1)$st order derivative
$$\frac{d^k g}{ds^k} + R_2(s) \frac{d^{k-2} g}{ds^{k-2}} + \ldots + R_k(s) g = 0 \ , \ R_i(s) \in \mathbb{C}(s)$$
whose fundamental solutions define the same projective period map as that of equation (\ref{eqnPF}).  
\end{definition}
It is always possible to pass to the projective normal form differential equation by rescaling each fundamental solution by the $k$th root of the Wronskian of the original equation.
\end{remark}

\begin{example}
Suppose now that $k = 2$, i.e., the initial differential equation is 
$$\frac{d^2 f}{ds^2} + P_1(s) \frac{df}{ds} + P_2(s) f = 0 \ ,$$
then the projective normal form of this differential equation takes the particularly simple form
$$\frac{d^2 g}{ds^2} + \left( P_2(s) - \frac{1}{2} {P_1}'(s) - \frac{1}{4} {P_1(s)}^2 \right) g = 0 \ .$$
\label{egPNF2}
\end{example}

As the process of taking the projective normal form does not alter the position or type of a maximal unipotent monodromy point, and as the projective solution determines the mirror map, we find
\begin{lemma}
The mirror map of a family of Calabi-Yau manifolds about a point of maximal unipotent monodromy of the Picard-Fuchs equation is determined by the projective normal form of this differential equation.
\end{lemma}

\section{Orbifold uniformization}

In this section we will recall the classical theory of {\em orbifold uniformization} of Riemann surfaces by Fuchsian second order ordinary differential equations, including some basic facts about the Schwarzian derivative, orbifold uniformizing differential equations, and automorphic functions.  More details of much of the material presented here can be found in the nice book by Yoshida \cite{Yos}.

\subsection{Uniformization of orbifolds}

\mbox{}

\begin{definition}
Let $X$ be a complex manifold, $Y \subset X$ a hypersurface, $Y = \cup_j Y_j$ its decomposition into irreducible components.  Let the numbers $b_j$, either $\infty$ or an integer $\geq 2$, be the {\em weights} attached to the corresponding $Y_j$.  The triple $(X, Y, \vec{b})$ is called an {\em orbifold} if for every point in $X \setminus \cup_j \{Y_j \ | \ b_j = \infty \}$ there is an open neighborhood $U$ and a covering manifold which ramifies along $U \cap Y$ with the given indices $\vec{b}$.  
\end{definition}

\begin{definition}
The orbifold $(X, Y, \vec{b})$ is called {\em uniformizable} if there is a simply connected global covering manifold, or {\em uniformization}, of $X$ with the given ramification data $(Y, \vec{b})$. 
\end{definition}

\begin{remark}
Often when one speaks of uniformization in classical settings it is understood that all the $b_j = \infty$ (i.e., only cusps occur).  We distinguish this from the definition given above by calling this more restrictive notion {\em classical uniformization}.  For emphasis we will often refer to uniformization where $b_j$ is allowed to be finite as well as $\infty$ as {\em orbifold uniformization}.
\end{remark}

\begin{definition}
Let $(X, Y, \vec{b})$ be an orbifold and $M$ its uniformization.  The multivalued inverse map $X \rightarrow M$ of the projection $M \rightarrow X$ is called the {\em developing map}, uniquely determined up to the group $Aut(M)$ of holomorphic automorphisms of $M$.
\end{definition}

For our application, we take $X = \mathbb{P}^1(\mathbb{C})$, $Y =$ a finite set of points, $M =$ the upper half plane $\mathbb{H}$, and $Aut(M) = PSL(2, \mathbb{R}) \simeq PU(1,1)$.

\subsection{Schwarzian derivatives, uniformizing differential equations, and automorphic functions}

\label{ssorbuni}

\mbox{}

\begin{definition}
The {\em Schwarzian derivative} $\{w;x\}$ of a nonconstant smooth function $w(x)$, with respect to $x$, is defined to be
$$\{w;x\}  =  \frac{3 (w'')^2 - 2 w' w'''}{4 (w')^2} \ .$$
\end{definition}

\begin{proposition}[\cite{Yos}, Proposition 4.1.1]
The Schwarzian derivative satisfies ``$PGL(2, {\mathbb{C}})$-invariance'':
$$\left\{ \frac{a w + b}{c w + d} \ ; \ x \right\} = \{w;x\} \ , \ \mbox{for all} \ \left( \begin{array}{cc} a & b \\ c & d \end{array} \right) \in PGL(2, {\mathbb{C}})$$
\label{propschwarz}
\end{proposition}

\begin{definition}
A function $w(x)$ is called {\em $PGL(2, {\mathbb{C}})$-multivalued} if any two branches of $w(x)$ are related by a projective transformaion.
\end{definition}

An immediate corollary of Proposition \ref{propschwarz} is
\begin{corollary}
If a function $w=w(x)$ is $PGL(2, {\mathbb{C}})$-multivalued, then the Schwarzian derivative $\{w;x\}$ is single valued.
\end{corollary}

\begin{proposition}[\cite{Yos}, Proposition 4.2]
Let $x$ and $w$ be the coordinates of $X$ and $M$ respectively.  
There are two linearly independent solutions $u_1$ and $u_2$ of the equation
\begin{eqnarray}
\frac{d^2u}{dx^2} - \{w(x);x\} u &=& 0
\label{eqnschude}
\end{eqnarray}
with single valued coefficients such that $w(x) = u_1(x) / u_2(x)$.
\label{propude}
\end{proposition}

\begin{definition}
The differential equation (\ref{eqnschude}) is called the {\em (orbifold) uniformizing differential equation}, or {\em Schwarzian differential equation}, of the orbifold $(X,Y,\vec{b})$.
\end{definition}

An immediate corollary of Proposition \ref{propude} is
\begin{corollary}
The projective solution of the uniformizing differential equation of the orbifold $(X,Y,\vec{b})$ is the developing map $w(x)$.
\end{corollary}

Let $\Gamma \subset PSL(2, {\mathbb{R}})$ be a discrete subgroup such that $\Gamma \, \backslash  \, {\mathbb{H}}^*$ is compact (i.e., a Fuchsian group of the first kind).  In our applications we will be primarily interested in a particular class of orbifold uniformizing differential equations associated to such groups $\Gamma$ with $\Gamma \, \backslash \, \mathbb{H}^*$ of genus zero.

\begin{definition}
A  function $f(\tau)$ on the upper half plane ${\mathbb{H}}$ is an {\em automorphic function} for $\Gamma$ if 
$$f(\gamma \tau) = f(\tau) \ , \ \mbox{for all} \ \gamma \in \Gamma \ \mbox{and} \ \mbox{for all} \ \tau \in {\mathbb{H}} \ ,$$
(i.e., invariance) and there are well defined limits of $f(\tau)$ as $\tau$ approaches parabolic vertices (i.e., $f(\tau)$ is meromorphic at each cusp of $\Gamma$).
\end{definition}

\begin{definition}
A univalent automorphic function $f$ for a genus zero Fuchsian group is called a {\em Hauptmodul}.
\end{definition}
The covering map for the uniformization $\mathbb{H} \rightarrow \Gamma \, \backslash \, \mathbb{H}^*$ is a Hauptmodul if $\Gamma$ is genus zero.

\begin{theorem}[\cite{Leh}, Theorem 6A]
Let $\Gamma$ be a Fuchsian group of the first kind.  If $f, g$ are automorphic functions for $\Gamma$, then they satisfy an algebraic equation
$$\mathcal{P}(f,g) = 0$$
with complex coefficients.
\end{theorem}

\begin{definition}
Such an algebraic equation satisfied by a pair of automorphic functions $f$ and $g$ for $\Gamma$ is known classically as a {\em modular equation}. Warning: Often in the literature this name is reserved for the particular case of modular equations relating the modular function $J(\tau)$ and $J(n \cdot \tau)$, $n \in \mathbb{N}$.
\end{definition}

\begin{corollary}[\cite{Leh}, Theorem 6B]
If there is a univalent automorphic function $f$ for $\Gamma$ (i.e., one to one on the fundamental domain), then every automorphic function for $\Gamma$ is a rational function of $f$.
\label{corlehner}
\end{corollary}

Suppose we are given a pair of genus zero groups $\Phi \subseteq \Gamma$ (finite index) with respective Hauptmoduln $f$, $g$.  Then the rational function $\mathcal{R}(z)$ from Corollary \ref{corlehner} expresses $g$ in terms of $f$: 
$$g = \mathcal{R}(f) \ .$$
\begin{definition}
The rational function $\mathcal{R}(z)$ is called the {\em modular relation} between the Hauptmoduln $f$ and $g$.
\end{definition}
For several examples of modular relations between $J(\tau)$ and Hauptmoduln of certain modular groups with level structure, see \cite{CY}.

\begin{remark}
The orbifold uniformizing differential equation is, up to projective equivalence, the ``pullback by $\mathcal{R}$'' of that for $\Gamma \, \backslash \, \mathbb{H}^*$.
In the next subsection we describe necessary and sufficient conditions on a rational function for it to pullback one uniformizing differential equation to (the projective equivalence class of) another.
\label{remarkitems}
\end{remark}

\subsection{Regular singular points and orbifold uniformization}

\label{regsingp}

\mbox{}

There is a particularly simple characterization of projective normal form equations which orbifold uniformize, given in terms of their characteristic exponents:
\begin{lemma}[\cite{Yos}, p.~50]
A second order Fuchsian differential equation in projective normal form is an orbifold uniformizing differential equation if and only if at each of its regular singular points the difference in characteristic exponents is either $0$ or $1/b$ for some $b \geq 2 \in \mathbb{Z}$.
\label{lemmayos50}
\end{lemma}

As we have already indicated, it is important that we understand the behavior of uniformizing differential equations under pullback.
Suppose then that we pullback a second order Fuchsian equation 
$L_2 f = 0$, with
$$L_2 = \frac{d^2}{dx^2} + P_1(x) \frac{d}{dx} + P_2(x) \ ,$$ 
by a rational function $x = \mathcal{R}(z)$.  

\begin{proposition}
The regular singular points of the pulled back equation $\mathcal{R}^*(L_2)$ lie among the inverse image points of the regular singular points of the original equation $L_2$ and the extra ramification points of the rational function $\mathcal{R}(z)$.
\end{proposition}
\begin{proof}
This is a straightforward computation:
$$\frac{df}{dx} = \frac{df}{dz} \left/ \frac{d\mathcal{R}}{dz} \right. $$
and
$$\frac{d^2f}{dx^2} = \left( \frac{d\mathcal{R}}{dz} \right)^{-2} \left( \frac{d^2f}{dz^2} - \frac{df}{dz} \left( \frac{d^2 \mathcal{R}}{dz^2} \left/ \frac{d \mathcal{R}}{dz} \right. \right) \right)$$
so that the pullback $\mathcal{R}^*(L_2) = 0$ of the equation $L_2 f = 0$ takes the form
$$\frac{d^2f}{dz^2} + \frac{df}{dz} \left( P_1(\mathcal{R}(z)) \frac{d \mathcal{R}}{dz} - \left( \frac{d^2 \mathcal{R}}{dz^2} \left/ \frac{d \mathcal{R}}{dz} \right. \right) \right) + f \left( P_2(\mathcal{R}(z)) \left(\frac{d \mathcal{R}}{dz} \right)^2 \right) = 0 \ .$$
Thus the only poles in the coefficients can occur at points where $P_1(\mathcal{R}(z))$ has one, or where $P_2(\mathcal{R}(z))$ has one, or at other places where $d \mathcal{R} / dz = 0$ (i.e., the {\em extra} ramification points of the $\mathcal{R}$ map).
\end{proof}

\begin{lemma}
The characteristic exponents of a point in the inverse image as above are exactly those of the point being pulled back times the order of ramification.  The characteristic exponents at an extra ramification point differ by an integer.
\label{lemmaramprod}
\end{lemma}
\begin{proof}
For simplicity assume that the singular point of $L_2$ under consideration is at $x = 0$.  In general (in the absence of ramification) there are $\deg \mathcal{R}(z)$ points in the inverse image of $0$, each with the same characteristic exponents as has $L_2$ at $0$.  If ramification of $\mathcal{R}(z)$ occurs at an inverse image point $z_0$, then $\mathcal{R}(z)$ has a multiple root there with multiplicity equal to the order of ramification.  Thus the characteristic exponents at $z_0$ are those of $L_2$ at $0$ times this order of ramification.  At an {\em extra} ramification point $z_1$, i.e., a point in the inverse image of an ordinary point of $L_2$ (characteristic exponents $0, 1$), the characteristic exponents are thus $0$ and the ramification index of $\mathcal{R}(z)$ at $z_1$.
\end{proof}

In light of Lemmas \ref{lemmayos50} and \ref{lemmaramprod}, we have 
\begin{theorem}
Characterization of types of regular singular points of a pullback of an orbifold uniformizing differential equation by a rational function:
\begin{itemize}
\item The extra ramification points of the rational function are all apparent singularities.
\item The points in the inverse image of logarithmic singularities of orbifold type (i.e., those with equal characteristic exponents) are logarithmic singularities of orbifold type.
\item For a regular singular point in the inverse image of a finite order orbifold point with characteristic exponent difference $1/b$:  
if the characteristic exponent difference is integral (i.e., if $b$ divides the ramification index), then the regular singular point is apparent;
if the characteristic exponent difference is not integral (i.e., if $b$ does not divide the ramification index), then the singular point is generic.
\end{itemize}
\label{thpullcharexp}
\end{theorem}

\begin{remark}
Some comments on the proof: The regular singular points which are extra ramification points have characteristic exponent difference an integer $\geq 2$ (if the difference were equal to $1$ this would be an ordinary point \cite[p.~23]{Yos}).  
The characteristic exponent difference is rescaled by the degree of the rational function; if this difference is zero it remains so.  
In the last case, we can say more about the characteristic exponent difference: by Lemma \ref{lemmaramprod} we know that this is the ramification degree $r$ divided by $b$.  The reduced form of this fraction has numerator $r / \gcd(r,b)$.  If this equals $1$, then such a point satisfies the orbifold uniformization criterion with finite weight $b/ \gcd(r,b)$ for the pullback equation.
\end{remark}

Theorem \ref{thpullcharexp} can be interpreted as providing a characterization of the modular relations between Hauptmoduln of a given genus zero Fuchsian group of the first kind and its finite index subgroups.

\section{Elliptic curve case: elliptic modular surfaces}

We are now ready to consider the problem of orbifold uniformization by Picard-Fuchs equations for families of elliptic curves --- the first case of {\em Picard-Fuchs uniformization}.  We will obtain a simple characterization of elliptic surfaces with modular mirror maps in terms of Kodaira's {\em functional invariant}.  This rational function assigns to each point in the base of an elliptic fibration the $J$-invariant of the corresponding elliptic curve fiber.  Moreover, the projective normal form of the pullback by the functional invariant of the uniformizing differential equation for the elliptic modular function $J(\tau)$ is the projective normal form of the Picard-Fuchs equation of the elliptic fibration.  Thus the mirror map is determined by the functional invariant, and using Theorem \ref{thpullcharexp} we obtain a criterion for modularity of the mirror map in terms of properties of the functional invariant alone.  Finally, we indicate how to use the modular relations between Hauptmoduln for genus zero subgroups of $PSL(2,\mathbb{Z})$ and $J(\tau)$ to provide explicit algebraic expressions for elliptic modular surfaces with specified modularity properties.

\subsection{Weierstrass elliptic surfaces and Kodaira theory}

\mbox{}

\begin{definition}
Let $p: X \rightarrow C$ be a flat proper map from a reduced irreducible ${\mathbb{C}}$-scheme $X$ to a complete, smooth curve $C$, such that every geometric fiber is an irreducible curve of arithmetic genus one, i.e., each fiber is one of
\begin{enumerate}
	\item an elliptic curve, or
	\item a rational curve with a node, or
	\item a rational curve with a cusp.
\end{enumerate}
The total space
$X$ is normal, and the generic fiber of $p$ is smooth.  Assume further that a section $s : C \rightarrow X$ is given, not passing through the nodes or cusps of the fibers.  Call the collection of such data $(p : X \rightarrow C, s)$ a {\em Weierstrass fibration} over $C$.  We may resolve the singularities of $X$ to obtain an elliptic surface with section $\overline{p} : \overline{X} \rightarrow C$, called the {\em induced elliptic surface}.  
\end{definition}

In fact there is a {\em canonical form} for such a Weierstrass fibration, exhibiting $X$ as a divisor in a ${\mathbb{P}}^2$-bundle over the base curve $C$.  
\begin{theorem}[\cite{Mir1}, Theorem (2.1)]
Let $\Sigma$ denote the given section of $p$, i.e., $\Sigma = s(C)$, a divisor on $X$ which is taken isomorphically onto $C$ by $p$.
Let $\mathcal{L} = p_*[\mathcal{O}_X(\Sigma) / \mathcal{ O}_X]$.  Suppose that the general fiber of $p$ is smooth.  Then $\mathcal{L}$ is invertible and $X$ is isomorphic to the closed subscheme of ${\mathbb{P}} = {\mathbb{P}}(\mathcal{L}^{\otimes 2} \oplus \mathcal{L}^{\otimes 3} \oplus \mathcal{ O}_Y)$ defined by
$$y^2 z = 4 x^3 - g_2 x z^2 - g_3 z^3 \ ,$$
where 
$$g_2 \in \Gamma(C, \mathcal{L}^{\otimes -4}) \ , \ g_3 \in \Gamma(C, \mathcal{L}^{\otimes -6}) \ ,$$
and $[x,y,z]$ is the global coordinate system of ${\mathbb{P}}$ relative to $(\mathcal{L}^{\otimes 2}, \mathcal{L}^{\otimes 3}, \mathcal{O}_{C})$.  Moreover the pair $(g_2, g_3)$ is unique up to isomorphism, and the {\em discriminant}
$$g_2^3 - 27 g_3^2 \in \Gamma(C, \mathcal{L}^{\otimes -12})$$
vanishes at a point $q \in C$ precisely when the fiber $X_q$ is singular.
\end{theorem}

\begin{definition}
Let $\mathcal{ J}$ denote the composition of the period morphism $\omega_1 / \omega_0 : C \rightarrow {\mathbb{H}}$ and the morphism $J : {\mathbb{H}} \rightarrow {\mathbb{P}}^1$ extending the classical modular function. Kodaira calls this the {\em functional invariant} of $\pi : X \rightarrow C$, i.e., $\mathcal{J} = J \circ \omega_1 / \omega_0$.
\end{definition}

Kodaira has classified the singular fiber types which can arise in Weierstrass fibered elliptic surfaces by describing these configurations of rational curves.  
\begin{theorem}[\cite{Kod}]
The singular fibers which appear in a smooth minimal elliptic surface fall into ``types'': $I_n$ ($n \geq 0$), $II$, $III$, $IV$, $I_n^*$ ($n \geq 0$), $IV^*$, $III^*$, and $II^*$.  Let $I_0$ denote a smooth elliptic fiber.  The fiber of type $I_1$ is a rational curve with a single node.  More generally, fibers of type $I_n$ consist of a $n$-cycle of intersecting rational curves for $n \geq 1$.  A fiber of type $II$ is just a rational curve with a single cusp.  Type $III$ fibers consist of two rational curves with a single point of tangency.  Fibers of type $IV$ consist of three rational components intersecting at a single point.  There are also fibers of types $I_n^*$, $n \geq 0$, $IV^*$, $III^*$, and $II^*$, whose dual intersection graphs, minus in each case a multiplicity one component, correspond to those graphs of Dynkin types $D_{n+4}$, $E_6$, $E_7$, and $E_8$ respectively.
\end{theorem}

We now recall how the Kodaira fiber types correlate with the ramification behavior of the $\mathcal{ J}$-map.

\begin{lemma}[\cite{Mir2}, Lemma IV.4.1]
Let $G = X_q$ be the fiber of $\pi$ over $q \in C$.
\begin{enumerate}
\item If $G$ has type $II$, $IV$, $IV^*$, or $II^*$, then $\mathcal{ J}(q) = 0$.  Conversely, suppose that $\mathcal{J}(q) = 0$.  Then
\begin{itemize}
	\item $G$ has type $I_0$ or $I_0^*$ if and only if $m(\mathcal{ J}) \equiv 0 \mod{3}$.
	\item $G$ has type $II$ or $IV^*$ if and only if $m(\mathcal{ J}) \equiv 1 \mod{3}$.
	\item $G$ has type $IV$ or $II^*$ if and only if $m(\mathcal{ J}) \equiv 2 \mod{3}$.
\end{itemize}
\item If $G$ has type $III$ or $III^*$, then $\mathcal{J}(q) = 1$.  Conversely, suppose that $\mathcal{J}(q) = 1$.  Then
\begin{itemize}
	\item $G$ has type $I_0$ or $I_0^*$ if and only if $m(\mathcal{ J}) \equiv 0 \mod{2}$.
	\item $G$ has type $III$ or $III^*$ if and only if $m(\mathcal{ J}) \equiv 1 \mod{2}$.
\end{itemize}
\item $G$ has type $I_n$ or $I_n^*$ with $n \geq 1$ if and only if $\mathcal{ J}$ has a pole at $q$ of order $n$.
\end{enumerate}
\label{lemmaGJlist}
\end{lemma}

\subsection{Picard-Fuchs equations of Weierstrass elliptic surfaces}

\label{PFWei}

\mbox{}

In this subsection we will describe the Picard-Fuchs differential equations for a Weierstrass elliptic curve fibration.  We first recall Griffiths' approach, which starts with a first order regular system from which the second order equation can be determined.  This is computationally convenient, but will not be so useful for our later purposes.  Ultimately, we will only care about properties of the projective normal form of the Picard-Fuchs differential equation.  For that reason we also present Stiller's approach by ``pulling back by $\mathcal{J}$'' which suffices to determine this.

\begin{theorem}[\cite{Sas}, p.~304]
Griffiths' approach to computing Picard-Fuchs yields
$$\frac{d}{dz} \left( \begin{array}{c} \eta_1 \\ \eta_2 \end{array} \right) = \left[ \begin{array}{cc} \frac{-1}{12} \frac{d \Delta / d z}{\Delta} & \frac{3 \delta}{2 \Delta} \\ \frac{-g_2 \delta}{8 \Delta} & \frac{1}{12} \frac{d \Delta / d z}{\Delta} \end{array} \right] \left( \begin{array}{c} \eta_1 \\ \eta_2 \end{array} \right)$$
where
$$\Delta(z) = g_2(z)^3 - 27 g_3(z)^2 \ ,$$
$$\delta(z) = 3 g_3(z) \frac{d g_2(z)}{dz} - 2 g_2(z) \frac{d g_3(z)}{dz} \ ,$$
and
$$\eta_1 = \int_\gamma \frac{dx}{y} \ , \ \eta_2 = \int_\gamma \frac{x dx}{y} \ .$$
\end{theorem}

Now we turn to Stiller's ``pulling back by the $\mathcal{ J}$-map'' approach.
Let $X$ be a complete smooth connected algebraic curve over ${\mathbb{C}}$ with function field $\mathbb{C}(X)$. 
Recall that the uniformizing differential equation for the group $PSL(2, \mathbb{Z})$ is 
$$ \Lambda f = \frac{d^2 f}{dx^2} + \frac{36 x^2 - 41 x + 32}{144 x^2 (x-1)^2} f = 0 \ ,$$
the projective normal form of the Picard-Fuchs equation of the family $\mathcal{E}$ of Example \ref{egellJ}.  Denote by $\Lambda_{\mathcal{J}}$ the projective normal form of the pullback equation $\mathcal{J}^*(\Lambda)$.
By an explicit computation one can show that
\begin{proposition}
Given a Weierstrass fibration with functional invariant $\mathcal{ J}$, the projective normal form of the Picard-Fuchs differential equation and the equation $\Lambda_\mathcal{ J}$ are identical.
\label{PropPFLJ}
\end{proposition}
\noindent
Actually, this can also be seen indirectly through Stiller's classification of {\em $K$-equations} \cite{Sti}.

We now can recast the results of Lemma \ref{lemmaGJlist} above as describing the singular fiber types associated to various regular singular points of the projective normal form of Picard-Fuchs $\Lambda_\mathcal{J}$, using Theorem \ref{thpullcharexp}.

\begin{proposition}
A regular singular point of $\Lambda_\mathcal{ J}$ corresponds to a fiber of type
\begin{itemize}
\item $I_n$ or $I_n^*$, $n \geq 1$ $\Leftrightarrow$ the characteristic exponents are equal
\item $I_0$ or $I_0^*$ $\Leftrightarrow$ the difference in characteristic exponents is a nonnegative integer
\item $II$ or $IV^*$ $\Leftrightarrow$ the difference in characteristic exponents is $1/3$ modulo integers
\item $III$ or $III^*$ $\Leftrightarrow$ the difference in characteristic exponents is $1/2$ modulo integers
\item $IV$ or $II^*$ $\Leftrightarrow$ the difference in characteristic exponents is $2/3$ modulo integers
\end{itemize}
\end{proposition}

\subsection{Characterization and construction of elliptic surfaces with modular mirror maps}

\mbox{}

The main result for elliptic curve families is now an immediate consequence of the criterion for orbifold uniformization (Theorem \ref{thpullcharexp}) in the elliptic surface setting:
\begin{theorem}
The projective normal form of the Picard-Fuchs equation of a Weierstrass elliptic surface will orbifold uniformize if and only if there are no apparent singularities and the functional invariant ${\mathcal{J}}(t)$ has zeros to orders $= 1$ or $\equiv 0 \bmod{3}$ and ones to orders $= 1$ or $\equiv 0 \bmod{2}$.
\label{thmellcase}
\end{theorem}
This already generalizes results of Stiller \cite{Sti} on classical uniformization by Picard-Fuchs differential equations, and Shioda \cite{Shi} on elliptic modular surfaces.  In addition, this also answers a question of Harnad and McKay regarding the geometric nature of modular solutions to generalized equations of Halphen type \cite{HM}.

Furthermore, by interpreting geometrically the absence of extra ramification of the functional invariant, we have this purely geometric characterization of modularity:
\begin{corollary}
The mirror map $z(q)$ of a Weierstrass elliptic surface over ${\mathbb{P}}^1$ is modular if and only if it is isolated in deformations which preserve the Kodaira singular fiber types, and it contains {\em no} singular fibers of types $IV$ or $II^*$.
\label{corellcase}
\end{corollary}

Applying this result to the list of rational elliptic surfaces over ${\mathbb{P}}^1$ we obtain:
\begin{theorem}
Here is the complete list of rational elliptic modular surfaces, i.e., elliptic surfaces over $\mathbb{P}^1$ with Euler characteristic $12$, non-constant functional invariant, and modular mirror maps, listed according to their singular fiber types:
$$\begin{array}{lll}
I_1, II, III^*  	&  I_2, I_4, III, III	&  I_2, I_3, I_4, III    \\
I_2, II, IV^*		&  I_3, I_3, III, III	&  I_1, I_1, I_8, II	\\
I_1, I_2, III^*		&  I_1, I_6, II, III	&  I_1, I_2, I_7, II	\\
I_3, III, III, III	&  I_2, I_5, II, III	&  I_1, I_4, I_5, II	\\
I_1, I_3, IV^* 		&  I_3, I_4, II, III	&  I_2, I_3, I_5, II	\\
I_4, II, III, III	&  I_1, I_7, II, II	&  I_1, I_1, I_1, I_9	\\
I_5, II, II, III	&  I_2, I_6, II, II	&  I_1, I_1, I_2, I_8	\\
I_1, I_1, I_4^* 	&  I_4, I_4, II, II	&  I_1, I_2, I_3, I_6	\\
I_2, I_2, I_2^* 	&  I_1, I_1, I_7, III	&  I_1, I_1, I_5, I_5	\\
I_6, II, II, II		&  I_1, I_2, I_6, III	&  I_2, I_2, I_4, I_4	\\
I_1, I_5, III, III	&  I_1, I_3, I_5, III	&  I_3, I_3, I_3, I_3	\\ 
\end{array}$$
\end{theorem}
\begin{proof}
Persson and Miranda \cite{Per, Mir3} have cataloged the singular fiber types of all rational elliptic surfaces over $\mathbb{C}$.  Miranda's list is particularly convenient to use for this application because it indicates as well the amount of excess ramification of the functional invariant associated to each combination of singular fibers.  Between the constraint on the singular fibers and lack of excess ramification of $\mathcal{J}$ (both by Theorem \ref{thmellcase}), we are able to determine the above list of $33$ types.  Explicit expressions for the mirror map in each case can be derived from the Weierstrass canonical form of each elliptic surface.  Since the number of singular fibers present is $3$ or $4$, explicit Weierstrass presentations are available in the literature already: see Schmickler-Hirzebruch for the $3$ singular fibers cases \cite{S-H}, and Herfurtner for the surfaces with $4$ singular fibers \cite{Her}.
\end{proof}

It seems appropriate to make some additional comments on the construction of explicit algebraic elliptic modular surfaces over ${\mathbb{P}}^1$ as an application of these results. 

Associated to any elliptic surface meeting the equivalent criteria of Theorem \ref{thmellcase} and Corollary \ref{corellcase}, there is a mirror map which is a Hauptmodul (un-normalized!) of a genus zero (finite index) subgroup of $PSL(2,{\mathbb{Z}})$.  If we denote this Hauptmodul by $H$, then we know that the functional invariant $\mathcal{J}$ expresses the modular function $J$ in terms of $H$, $J(q) = \mathcal{J}(H(q))$, and the functional invariant defines a modular relation between the two genus zero modular groups.

We know that the uniformizing differential equation $\Lambda$ for $PSL(2,\mathbb{Z})$ is also the projective normal form of the Picard-Fuchs equation of the family of elliptic curves $\mathcal{E}$ in Example \ref{egellJ}.  We can identify the base parameter $s$ of that family with the modular function $J$.  The base parameter of the ``pullback surface'' obtained by replacing $J$ everywhere with $\mathcal{J}(H)$ in the equation for $\mathcal{E}$ is likewise identified with $H$.  Since we are provided  with an explicit algebraic elliptic modular surface with mirror map a Hauptmodul for $PSL(2, \mathbb{Z})$, it is a simple matter to determine an explicit algebraic elliptic surface with mirror map a Hauptmodul $H$ for any specified modular subgroup.  We will return to this sort of ``inverse construction'' at the end of the next section in the context of K3 surfaces and a question of Dolgachev.

\section{$M_n$-polarized K3 surface case: Mirror-Moonshine}

\label{Mnpol}

The next case of Picard-Fuchs uniformization is by one parameter families of certain {\em lattice polarized} K3 surfaces.  The direct motivation for considering this case is of course the Mirror-Moonshine Conjecture of Lian-Yau \cite{LY1}.  In the first subsection we recall the definition of $M_n$-polarized K3 surfaces and Dolgachev's results on their moduli.  Next we show how classical Fuchsian ordinary differential equation facts imply that the Picard-Fuchs equation of a family of such surfaces is a ``symmetric square''.  After a brief overview of the Mirror-Moonshine Conjecture, we carefully formulate the questions we want to answer in \S \ref{MMph}.  A crucial tool, the generalized functional invariant, is introduced in the next subsection.  Finally, in \S \ref{commenmodres} we prove our results on commensurability and modularity of the mirror map in the $M_n$-polarized K3 surface case.

\subsection{$M_n$-polarized K3 surfaces.}

\mbox{}

See the paper in this volume by Verrill and Yui \cite{VY} for a general discussion of lattice polarized K3 surfaces after Dolgachev \cite{Dol}.

Let $L$ be the K3 lattice
$$L = U \perp U \perp U \perp - E_8 \perp - E_8 \ .$$
Let $M = \langle 2n \rangle$.  Since the lattice $U$ represents any even integer, we can assume $M \subset U$.  Then 
$$M^\perp \simeq U \perp U \perp - E_8 \perp - E_8 \perp \langle - 2n \rangle =: M_n \ .$$
This lattice is denoted by $\mathcal{M}$ in \cite{VY}.

The global monodromy group of the ``universal family'' of K3 surfaces over the (coarse) moduli space of K3 surfaces polarized by a rank 19 lattice $M_n$ has been determined by Dolgachev.  Define to begin with
$$\Gamma_0(n) := \left\{ \left. \left( \begin{array}{cc} a & b \\ c & d \end{array} \right)  \in PSL(2, {\mathbb{Z}}) \ \right| c \equiv 0 \bmod n  \right\}  \subset PSL(2, {\mathbb{Z}})  \ .$$
\begin{definition} 
The element of order two
$$F_n = \left( \begin{array}{cc} 0 & -1 / \sqrt{n} \\ \sqrt{n} & 0 \end{array} \right) \in PSL(2, {\mathbb{R}})$$
in the normalizer of $\Gamma_0(n)$ in $PSL(2, {\mathbb{R}})$ is the {\em Fricke involution}.
\end{definition}

\begin{definition}
The {\em Fricke modular group} of level $n$ is the subgroup of $PSL(2, {\mathbb{R}})$ is generated by $\Gamma_0(n)$ and $F_n$, and denoted by $\Gamma_0(n)^*$.
\end{definition}

\begin{theorem}[\cite{Dol}, Theorem (7.1)]
The coarse moduli space of K3 surfaces polarized by $M_n$ is:
$$K_{M_n} \simeq {\mathbb{H}} / \Gamma_0(n)^* \ .$$
\label{ThDolMn}
\end{theorem}

For later convenience we recall
\begin{proposition}[\cite{Dol}, Theorem (7.5)]
Suppose that $X_0(n)^*$ is rational.  Then there exists a unique holomorphic function
$$H_n : {\mathbb{H}} \rightarrow {\mathbb{C}}$$
such that
\begin{enumerate}
\item $H_n$ is invariant with respect to $\Gamma_0(n)^*$
\item $H_n$ has a Fourier expansion
$$H_n(\tau) = \frac{1}{q} + \sum_{k=1}^\infty c_k q^k \ , \ q = e^{2 \pi \imath \tau}$$
\item the coefficients of the Fourier expansion are all integers
\item as a meromorphic function on $X_0(n)^*$, the function $h_n$ has a simple point at the cusp $\Gamma_0(n)^* \cdot \infty$ and generates the field of meromorphic functions on $X_0(n)^*$.
\end{enumerate}
\end{proposition}

\subsection{Picard-Fuchs equations of $M_n$-polarized K3 surface families}

\label{PFMn}

\mbox{}

Here we describe in some detail the symmetric square structure of the Picard-Fuchs equations of $M_n$-polarized K3 surface families.

\begin{theorem}[\cite{Sin}, Lemma 3.1.(b)]
Let $L_1(y)$ and $L_2(y)$ be homogeneous linear differential polynomials with coefficients in ${\mathbb{C}}(t)$.  Then there exists a homogeneous linear differential equation $L_3(y) = 0$ with coefficients in $\mathbb{C}(t)$ and solution space the $\mathbb{C}$-span of 
$$\{\nu_1 \nu_2 \ | \ L_1(\nu_1) = 0 \ \mbox{and} \ L_2(\nu_2) = 0 \} \ .$$
\end{theorem}

\begin{definition}
We call the operator $L_3(y)$ constructed above the {\em symmetric product} of $L_1$ and $L_2$, and denote it by $L_1 \circledS L_2$.  In fact, the operation is associative, and we may further define $L^{\circledS n}$ for $n \geq 1$ by $L^{\circledS 1} = L$  and $L^{\circledS n} = L^{\circledS n-1} \circledS L$.  We call $\mbox{Sym}^n(L) = L^{\circledS n}$ the {\em $n$th symmetric power} of $L$; conversely, $L$ is the {\em $n$th root} of $L^{\circledS n}$.
\end{definition}

\begin{lemma}[\cite{Sin}, Lemma 4, p.~129]
Let $L(y)$ be a homogeneous linear differential polynomial with coefficients in ${\mathbb{C}}(t)$.  Then $L(y) = L_2^{\circledS n}(y)$ for some second order homogeneous linear differential polynomial $L_2(y)$ with coefficients in $\mathbb{C}(t)$ if and only if there exists a fundamental set of solutions $\{y_1, \ldots, y_{n+1}\}$ of $L(y) = 0$ such that
$$y_i y_{i+2} - y_{i+1}^2 = 0 \ , \ i = 1, \ldots, n-1 \ .$$
\label{lemmasinger}
\end{lemma}

\begin{corollary}
Let $L(y) = 0$ be a third order homogeneous linear equation with coefficients in ${\mathbb{C}}(t)$.  If there exists a nondegenerate homogeneous polynomial $P$ of degree $2$ with constant coefficients  and a fundamental set of solutions $\{y_1, y_2, y_3\}$ of $L(y) = 0$ such that $P(y_1, y_2, y_3) = 0$, then $L(y)$ is the second symmetric power of a second order homogeneous linear differential equation with coefficients in ${\mathbb{C}}(t)$.
\label{corfano}
\end{corollary}
\begin{proof}
This follows easily from Lemma \ref{lemmasinger}.  By assumption, the fundamental set of solutions satisfies a nondegenerate quadratic relation.  Since all such quadrics in ${\mathbb{P}}^2({\mathbb{C}})$ are projectively equivalent to
$$y_1 y_3 - y_2^2 = 0$$
the criterion of the lemma applies and $L(y)$ is a symmetric square.
\end{proof}

Assume that we have a second order Fuchsian ordinary differential equation $L_2 f = 0$ where
$$L_2 = \frac{d^2}{dt^2} + P_1(t) \frac{d}{dt} + P_2(t) \ .$$
The second order equation $L_2 f = 0$ is equivalent to the system of first order differential equations
$$\left\{ \begin{array}{rcl}
f' & = & g \\
g' & = & - P_2 f - P_1 g 
\end{array} \right.$$
with $\{ f , g \}$ as fundamental solutions.  Observe that
$$\{ f^n , f^{n-1} g , \ldots , f g^{n-1} , g^n \}$$
forms a set of fundamental solutions for the $n$-th symmetric power $L = L_2^{\circledS n}$.  The following result describes a system of first order differential equations for $L$ with these fundamental solutions.
\begin{theorem}[\cite{Lee}, Theorem 2]
If $f$ and $g$ satisfy
$$\frac{d}{dt} \left( \begin{array}{c} f \\ g \end{array} \right) = \left( \begin{array}{cc} 0 & 1 \\ - P_2 & - P_1 \end{array} \right) \left( \begin{array}{c} f \\ g \end{array} \right) \ ,$$
then we have
$$\frac{d}{dt} \left( \begin{array}{c} f^n \\ f^{n-1} g \\ \vdots \\ f g^{n-1} \\ g^n \end{array} \right) = A \left( \begin{array}{c} f^n \\ f^{n-1} g \\ \vdots \\ f g^{n-1} \\ g^n \end{array} \right) \ ,$$
where $A = (a_{ij})$ is an $(n+1) \times (n+1)$ matrix such that
$$\begin{array}{rclcl}
a_{k,k} &=& (1-k) P_1 & \ , \ & 1 \leq k \leq n+1 \ , \\
a_{k, k+1} &=& n + 1 - k & \ , \ & 1 \leq k \leq n \ , \\
a_{k+1, k} &=& - k P_2 & \ , \ & 1 \leq k \leq n \ , \\
a_{i,j} &=& 0 & \ , \ & i > j + 1 \ \mbox{or} \ j > i + 1 \ .
\end{array}$$
\end{theorem}

\begin{example}
In particular, when $n = 2$, the case for a symmetric square, it one rewrites the system as a single third order equation
$$\mbox{Sym}^2(L_2) f = f''' + 3 P_1 f'' + (2 {P_1}^2 + 4 P_2 + {P_1}') f' + (4 P_1 P_2 + 2 {P_2}') f = 0 \ .$$
\end{example}

In this form it is easy to check, using the closed form expression for the projective normal form of a second order Fuchsian differential equation given in Example \ref{egPNF2}, that
\begin{proposition}
Let $L_2$ be as above a second order Fuchsian ordinary differential operator, and let $L = \mbox{Sym}^2(L_2)$ be its symmetric square.  Then the projective normal form of $L$ is the symmetric square of the projective normal form of $L_2$.
\end{proposition}

\begin{remark}
In fact, it is possible to provide an explicit description of the relationship between the monodromy matrices of the second order ``square root'' equation and those of the third order symmetric square equation.  This is provided by the  faithful representation of $SL(2, {\mathbb{C}})$ in $SL(3, {\mathbb{C}})$ via the symmetric square representation \cite{Sin}.
\end{remark}

Finally, we see the relevance of all of this for Picard-Fuchs equations of our $M_n$-polarized K3 surface families:
\begin{theorem}
The Picard-Fuchs equation of a family of $M_n$-polarized K3 surfaces is the symmetric square of a second order homogeneous linear Fuchsian ordinary differential equation.
\label{ThSym2K3}
\end{theorem}
\begin{proof}
To begin with, the order of the Picard-Fuchs equation is equal to the rank of the transcendental lattice, i.e., $22 - 19 = 3$.  By Dolgachev's Torelli theorem for lattice polarized K3 surfaces (specifically, the proof of Theorem \ref{ThSym2K3} above) the period domain lies on a nondegenerate quadric in ${\mathbb{P}}^2$.  Thus, Corollary \ref{corfano} implies that the third order Picard-Fuchs differential equation is in fact a symmetric square.
\end{proof}
There is another approach to proving this result which takes advantage of the special geometric properties of $M_n$-polarized K3 surfaces, namely their presentation as Shioda-Inose surfaces coming from a product of two elliptic curves linked by an $n$-isogeny.  See \cite{Pet, VY} for more details.

\subsection{The Mirror-Moonshine phenomenon}

\label{MMph}

\mbox{}

In their first systematic investigations of mirror symmetry for one parameter families of Calabi-Yau manifolds constructed via the ``orbifold construction'' \cite{LY1}, Lian and Yau discovered that the reciprocal of the mirror maps for the K3 surfaces they were studying agreed, up to an additive constant, with Thompson series in the lists of Conway-Norton \cite{CN}.  The evidence was sufficiently strong that they formulated
\begin{conjecture}[Lian-Yau, \cite{LY1, ICMP}]
If $z(q)$ is the mirror map for a one parameter deformation of an algebraic K3 surface from an orbifold construction which has a third order Picard-Fuchs equation, then, for some $c \in {\mathbb{Z}}$, the $q$-series
$$\frac{1}{z(q)} + c$$
is a Thompson series $T_g(q)$ for some element $g$ in the Monster.
\label{LYConj}
\end{conjecture}

In \cite{LY2}, Lian and Yau compute many more toric examples (including over a dozen complete intersection examples), and note that the correspondence to monstrous groups persists.  This suggested that the hypothesis regarding the ``orbifold construction'' should perhaps be weakened to the hypothesis ``torically constructed''.

As noted in the proof of Theorem \ref{ThSym2K3}, for a family of lattice polarized K3 surfaces the condition of having a third order Picard-Fuchs equation is equivalent to the family possessing a polarization by a lattice of rank 19.  In other words, the generic member of such a family will have Picard rank 19.  

Furthermore, recall that a Thompson series is in particular a Hauptmodul for some ``monstrous'' genus zero arithmetic group $\Gamma$, and that the various equivalent Hauptmoduln are well-defined as generators of the function field of the rational curve $\Gamma \, \backslash \,  \mathbb{H}^*$ only up to action of $\Gamma$.  We see that in Conjecture \ref{LYConj} an equivalent conclusion is that the mirror map is itself a Hauptmodul (unnormalized!) for some monstrous $\Gamma$.

Before Conjecture \ref{LYConj} was even formulated, Beukers, Peters, and Stienstra had computed the Picard-Fuchs equation of a particular family of $M_n$-polarized K3 surfaces \cite{BP, PS}.  The mirror map was determined by Verrill and Yui  \cite[(V.2)]{VY}.  It provides a counterexample to a ``monstrous'' generalization of the mirror-moonshine conjecture for torically constructed cases; although the mirror map is {\em commensurable} to a monstrous Hauptmodul, it is not a monstrous Hauptmodul itself.  

This suggests the following two questions:
\begin{description}
\item[Commensurability] {\em When is the mirror map commensurable to a Thompson series?}

\item[Modularity] {\em Can we characterize the families of $M_n$-polarized K3 surfaces whose mirror maps are Hauptmoduln for genus zero subgroups of $\Gamma_0(n)^*$?}
\end{description}
We will provide answers to each of these questions in \S \ref{commenmodres}.

\subsection{Generalized functional invariant}

\label{genfct}

\mbox{}

\begin{definition}
By analogy with the functional invariant of an elliptic surface, to any family of $M_n$-polarized K3 surfaces over a curve $C$ we associate the composition of the period morphism
$$\begin{array}{rcl} 
\Omega : C & \rightarrow & {\mathbb{H}} \\
z & \mapsto & \Omega_1(z) / \Omega_0(z)
\end{array}$$
and the canonical projection onto the coarse moduli space
$$ H_n : {\mathbb{H}}  \rightarrow  {\mathbb{H}} / \Gamma_0(n)^* \ ,$$
and call it the {\em generalized functional invariant} $\mathcal{H}_n = H_n \circ \Omega$.
\end{definition}

\begin{lemma}
The generalized functional invariant $\mathcal{H}_n$ is meromorphic.
\end{lemma}
\begin{proof}
This follows directly from the regularity of the singular points of Fuchsian differential equations and the meromorphicity of the Hauptmodul.  The question is local; we will always assume we work with charts centered about $0$.  Consider the second order (projective normal form of square root of) Picard-Fuchs differential equation on ${\mathbb{P}}^1$ with coordinate $z$.  Let $a_\rho$ denote the regular singular points.  Let $\Omega(z)$ be the multivalued holomorphic function from the complement of the $a_\rho$ to the upper half plane, $H_n$ be the Hauptmodul viewed as a function from $\mathbb{H}$ to $\mathbb{P}^1$, and $\mathcal{H}_n$ be the generalized functional invariant.  There are two cases to consider: the points $a_\rho$ of finite or of infinite order monodromy.  

If a point has finite order $m$ monodromy, then as a function of $z = \sigma^m$ the period $\Omega(\sigma^m)$ is a single-valued holomorphic function of $\sigma$ in a neighborhood $0 < \sigma < \epsilon$, where $\Im ( \Omega(\sigma^m) ) > 0$.  In particular by Picard's theorem the singularity is not essential. 

If the monodromy about an $a_\rho$ is of infinite order, we can argue that the singularity is not essential by bounding the rate of growth of $\mathcal{H}_n(z) = H_n(\Omega(z))$.  For appropriately chosen coordinate $w$ centered about the singularity, by Frobenius' method we know that $w = O(z^{-m} \log z)$, so $|w| < O(z^{-(m+1)})$ as $z \rightarrow 0$.  Because $H_n(w)$ is a meromorphic function on $\Gamma_0(n)^* \, \backslash \, {\mathbb{H}}^*$ (and so on ${\mathbb{H}}^*$), near a singular point of $H_n$ the rate of growth of $H_n(w)$ equals $O(w^{-r})$ for some finite $r > 0$.  Thus $|H_n(\Omega(z))| < O(z^{-r(m+1)})$, and ${\mathcal{H}}_n(z)$ does not have an essential singularity.  

Thus the singularity is either removable or a pole, and ${\mathcal{H}}(z)$ is meromorphic.
\end{proof}

\subsection{Commensurability and modularity results}
\label{commenmodres}

\mbox{}

We can now address the two main questions regarding the Mirror-Moonshine conjecture for any family of $M_n$-polarized K3 surfaces over ${\mathbb{P}}^1$:  The first result explains the relationship of the mirror map to Thompson series.  The second characterizes the modularity properties of the mirror map.

\subsubsection{Commensurability.}

\mbox{}

By checking the tables of Conway-Norton \cite{CN}, we find:  If $\Gamma_0(n)^*$ is a genus zero group with $n \neq 49, 50$, then it is monstrous, i.e., its Hauptmodul is a Thompson series.

Suppose we are given an algebraic family of $M_n$-polarized K3 surfaces over ${\mathbb{P}}^1$.  In particular the genus of $\Gamma_0(n)^*$ must be zero.  The generalized functional invariant is a rational function, $\mathcal{H}_n(z)$, such that
$$H_n(q) = \mathcal{H}_n(z(q)) \ ,$$
i.e., it expresses the Hauptmodul for $\Gamma_0(n)^*$ as a rational function of the mirror map $z(q)$.  If $n \neq 49, 50$, this means that there is a rational expression for a Thompson series $H_n(q)$ in terms of the mirror map $z(q)$.  In particular this shows that: 
\begin{lemma}
The mirror map of a family of $M_n$-polarized K3 surfaces over $\mathbb{P}^1$ is commensurable with a Thompson series when $n \neq 49, 50$.
\end{lemma}

The groups $\Gamma_0(49)^*$ and $\Gamma_0(50)^*$ are two of Conway-Norton's ``ghosts'', i.e., non-monstrous genus zero subgroups of $\Gamma_0(n)^*$ containing $\Gamma_0(n)$.  Their Hauptmoduln are however still ``replicable functions'' in the sense of Conway-Norton, and lie among the list of ``monster-like'' groups of Ferenbaugh \cite{Fer}.  So we can say:  
\begin{lemma}
The mirror map of any family of $M_n$-polarized K3 surfaces over $\mathbb{P}^1$ is commensurable with a replicable function (i.e., a Hauptmodul of a monster-like group).
\end{lemma}

\subsubsection{Modularity.}

\mbox{}

We wish now, by analogy with the elliptic curve case, to apply the general criterion for  pulling back orbifold uniformizing differential equations (Theorem \ref{thpullcharexp}) to $M_n$-polarized K3 surface families.  We replace the functional invariant with the rational function $\mathcal{H}_n(z)$, our generalized functional invariant.  The differential equation $\Lambda_{\mathcal{J}}$ is replaced by the projective normal form of the square root of the Picard-Fuchs differential equation of our one parameter family of $M_n$-polarized K3 surfaces.

As the uniformizing group for the ``target'' is now $\Gamma_0(n)^*$, what order elliptic points must be considered?  In \cite{Fer} it is indicated (page 55) that for the broader class of groups $G = \Gamma_0(n|h)+ e_1, e_2, \ldots$, an elliptic point can only have order $2$, $3$, $4$, or $6$.  We need to check that the ``order of value'' $p$ of $\mathcal{H}_n(z)$ (i.e., the vanishing order of $\mathcal{H}_n(z) - p$) for each elliptic point satisfy the conditions of Theorem \ref{thpullcharexp}.

\begin{theorem}
The square root of the projective normal form of the Picard-Fuchs equation of a one parameter family of $M_n$-polarized K3 surfaces (with genus zero $\Gamma_0(n)^*$) will orbifold uniformize a Zariski open subset of the base of the family if and only if it has no apparent singularities, and for each elliptic point $p \in X_0(n)^*$, $\mathcal{H}_n(z) - p$ vanishes to permissible orders as listed in Table \ref{tableMn}.
\label{LYchar}
\end{theorem}

\begin{table}[h]
\begin{center}
\begin{tabular}{| c | l |}
\hline
Order of elliptic point $p \in X_0(n)^*$ & Vanishing order $\mathcal{H}_n(z) - p$ \\ \hline
$2$ & $\equiv 0 \bmod 2$ or $= 1 \bmod 2$ \\
$3$ & $\equiv 0 \bmod 3$ or $= 1 \bmod 3$ \\
$4$ & $\equiv 0 \bmod 4$ or $= 1$ or $2 \bmod 4$ \\
$6$ & $\equiv 0 \bmod 6$ or $= 1$, $2$, or $3 \bmod 6$ \\ \hline
\end{tabular}
\end{center} 
\caption{Vanishing orders permitting modular mirror maps}
\label{tableMn}
\end{table}

\begin{proof}
The argument here is virtually identical to that for elliptic curves.  To begin with we work with the square root of the projective normal form of Picard-Fuchs, a second order Fuchsian ordinary differential equation. Replace $J(\tau)$ with the normalized hauptmodul $H_n(\tau)$ for the appropriate $\Gamma_0(n)^*$.  Replace the functional invariant $\mathcal{J}$ with the generalized functional invariant $\mathcal{H}_n$.  The inverse image of a cusp under the generalized functional invariant is still a cusp and these are all the cusps (finite index condition).  A regular singular point in the inverse image of an elliptic point is either an algebraic singularity, or an apparent singularity.  The algebraic singularities are permitted if and only if they are of the types listed in Table \ref{tableMn}.  As before, the ``extra'' points at which $d \mathcal{H}_n / dz = 0$ are just the apparent singularities outside of the inverse image locus, or equivalently those for which the characteristic exponents are distinct.  Positivity and surjectivity of the period morphism, and the fact that the universal cover of the complement of the (non-apparent) singularities is ${\mathbb{H}}$ all follow from the same properties of the $X_0(n)^*$ orbifold.
\end{proof}

We have thus characterized the Mirror-Moonshine phenomenon completely.  This characterization has the advantage that it is expressed purely in terms of properties of the generalized functional invariant, and applies equally well to any example family, toric or not.

\begin{remark}
As an application we can construct families of $M_n$-polarized K3 surfaces which orbifold uniformize with specified monodromy groups.  The only constraint is that we must already have at hand an algebraic construction for a family of $M_n$-polarized K3 surfaces with mirror map a Hauptmodul for the full group $\Gamma_0(n)^*$.  This thus provides a partial answer to a question of Dolgachev regarding the ``inverse problem'' (see \cite[(VI.2)]{VY}), reducing it to the case $G = \Gamma_0(n)^*$.  In other words, the problem can be solved completely just by finding explicit algebraic constructions of families of $M_n$-polarized K3 surfaces with mirror maps equal to $H_n(q)$.

As in the case of elliptic curve families, the starting point is that we know how to pullback from the ``universal family'' with a given lattice polarization, to another family with the same lattice polarization which orbifold uniformizes with a proper subgroup as monodromy group.  Assume for now that the subgroup is also of genus zero.

As long as we have an algebraic presentation (e.g., as hypersurfaces, complete intersections) for the universal family in the first place, then we can construct the pullback family.  The necessary pullback function is the rational function which defines the modular relation for the Hauptmodul of the universal family as a rational function in the Hauptmodul of the second family (i.e., for the finite index subgroup).  Once again the general fact (Corollary \ref{corlehner} and subsequent remarks) relating automorphic functions for genus zero groups provides this rational function.
\end{remark}

\begin{remark}
A necessary condition on the projective normal form Picard-Fuchs equation for orbifold uniformization once again is that there be {\em no apparent singularities}.  There is a great deal of evidence to support the conjecture that no apparent singularities arise in any torically constructed family of Calabi-Yau manifolds.  Aside from the lack of counterexamples, there are other reasons to believe such a conjecture (see the author's work on {\em geometric isomonodromic deformations}, \cite{Dor}).
\end{remark}

\section{Generalizations}

\label{general}

We comment now on generalizations of these Picard-Fuchs uniformization results in various directions, involving changes in the natures of both the fiber and the base of the families considered.  In each case, a complete  discussion may be found in \cite{Dor2}.

To begin with, we note that our results for $M_n$-polarized K3 surface families immediately extend to the general rank 19 lattice polarized case.  In particular, the arguments from Fuchsian theory for the symmetric square structure of the Picard-Fuchs equation are independent of the particular rank 19 lattice chosen.  What is lacking in the general setting is a complete classification of the rank 19 lattices and a specific description of monodromy groups.

It is possible to lift the convenient {\em genus zero} restriction on both the uniformizing groups and the bases of the families under consideration.  Most of the arguments given in determining the Picard-Fuchs uniformization criteria apply immediately to the case of families of elliptic curves and $M_n$-polarized K3 surfaces over higher genus curves.  These criteria still imply that the mirror map is an automorphic function for a (no longer genus zero) subgroup of the genus zero groups $PSL(2, \mathbb{Z})$ or $\Gamma_0(n)^*$ respectively, even though it is no longer a Hauptmodul.  Modification of the arguments to allow for $\Gamma_0(n)^*$ of genus higher than zero can be made as well.

Furthermore, there is an extension of our modularity criteria to certain multi-parameter families of abelian varieties and K3 surfaces of various Picard ranks.  In particular, this can be applied to two-parameter families of K3 surfaces with Picard rank 18 as studied in \cite[(VI.4)]{VY}.  This case is particularly nice, as it is possible to apply another old result of Fano to conclude that generic one-parameter ``slices'' of such families will possess Picard-Fuchs equations of order four, built from two distinct second order equations by {\em tensor product}.  When the two-parameter family is ``modular'', with respect to the appropriate Hilbert modular group, the Picard-Fuchs equations on the one-parameter slices reflect this by being ``bi-modular'', i.e., each second order factor equation is independently modular.  The example of Hosono, \cite[(VI.4.1)]{VY}, provides a nice illustration of this phenomenon.

For more on all of this, together with a discussion of the considerably more subtle case of Calabi-Yau threefold families, see \cite{Dor2}.


\bibliographystyle{amsalpha}

\begin{thebibliography}{ABC}

\bibitem [AB]{AB} D. Anosov and A. Bolibruch, \textit{The Riemann-Hilbert Problem}, Aspects of Mathematics, Friedr. Vieweg \& Sohn (1994). 

\bibitem [BPV]{BPV} W. Barth,  C. Peters   and  A.  Van de Ven, \textit{Compact Complex Surfaces}, Springer-Verlag (1984).

\bibitem [BP]{BP}  J. Bertin and C. Peters, \textit{Variations of Hodge Structures, Calabi-Yau Manifolds, and Mirror Symmetry}, Introduction \`{a} la Th\'{e}orie de Hodge, Panoramas et Synth\`{e}ses, 3, Paris: Soci\'{e}t\'{e} Math\'{e}matique de France (1996).

\bibitem [CY]{CY} I. Chen and N. Yui, \textit{Singular Values of Thompson Series}, Groups, Difference Sets, and the Monster: Proceedings of a Special Research Quarter at the Ohio State University, Spring 1993, Ohio State University Mathematical Research Institute Publications 4, Walter de Gruyter, New York (1996).

\bibitem [CK]{CK}   D. Cox and S. Katz, \textit{Mirror Symmetry and Algebraic Geometry}, Mathematical Surveys and Monographs, American Mathematical Society (1999).

\bibitem [CN]{CN} J. Conway and S. Norton, \textit{Monstrous Moonshine}, Bull. London Math. Soc. 11 (1979), no. 3, 308--339.

\bibitem [Dol]{Dol}  I. Dolgachev, \textit{Mirror Symmetry for Lattice Polarized $K3$ Surfaces}, Algebraic Geometry, 4. J. Math. Sci. 81 (1996), no. 3, 2599--2630.

\bibitem [Dor]{Dor} C. F. Doran, \textit{Picard-Fuchs Uniformization and Geometric Isomonodromic Deformations: Modularity and Variation of the Mirror Map}, Dissertation, Harvard University, Spring 1999.

\bibitem [Dor2]{Dor2} C. F. Doran, \textit{Picard-Fuchs Uniformization and Modularity of the Mirror Map}, in preparation.

\bibitem [EMM]{EMM}  \textit{Essays on Mirror Manifolds}, Ed. by S.-T. Yau, International Press, Hong Kong (1992).

\bibitem [Fer]{Fer}  C. Ferenbaugh, \textit{The Genus-Zero Problem for $n\vert h$-Type Groups}, Duke Math. J. 72 (1993), no. 1, 31--63.

\bibitem [HM]{HM} J. Harnad and J. McKay, \textit{Modular Solutions to Equations of Generalized Halphen Type}, solv-int/9804006.

\bibitem [Her]{Her}  S. Herfurtner, \textit{Elliptic Surfaces with Four Singular Fibres}, Math. Ann. 291 (1991), no. 2, 319--342.

\bibitem [HLY]{HLY} S. Hosono, B. Lian and S.-T. Yau, \textit{Maximal Degeneracy Points of GKZ Systems}, J. Amer. Math. Soc. 10 (1997), no. 2, 427--443. 

\bibitem [ICMP]{ICMP}  B. Lian and S.-T. Yau, \textit{Mirror symmetry, Rational Curves on Algebraic Manifolds and Hypergeometric Series}, Proceedings of International Congress of Mathematical Physics (1994).

\bibitem [Kod]{Kod}  K. Kodaira, \textit{On Compact Analytic Surfaces,  II}, Ann. of Math. (2) 77 (1963), 563--626.

\bibitem [Lee]{Lee}  M.-H. Lee, \textit{Picard-Fuchs Equations for Elliptic Modular Varieties}, Appl. Math. Letters 4 (1991), no. 5, 91--95.

\bibitem [Leh]{Leh}   J. Lehner, \textit{Discontinuous Groups and Automorphic Functions}, Mathematical Surveys, Number VIII, American Mathematical Society (1964).

\bibitem [LY1]{LY1} B. Lian and S.-T. Yau, \textit{Arithmetic Properties of Mirror Map and Quantum Coupling}, Comm. Math. Phys. 176 (1996), no. 1, 163--191.

\bibitem [LY2]{LY2}  B. Lian and S.-T. Yau, \textit{Mirror Maps, Modular Relations and Hypergeometric Series II}, Nucl. Phys. Proc. Suppl. 46 (1996) 248-262.

\bibitem [Mir1]{Mir1} R. Miranda, \textit{The Moduli of Weierstrass Fibrations over ${\mathbb{P}}^1$}, Math. Ann., 255, no. 3, 379--394 (1981).


\bibitem [Mir2]{Mir2} R. Miranda, \textit{The Basic Theory of Elliptic Surfaces}, Dottorato di Ricerca in Matematica. ETS Editrice, Pisa (1989).


\bibitem [Mir3]{Mir3} R. Miranda, \textit{Persson's List of Singular Fibers for a Rational Elliptic Surface}, Math. Z. 205 (1990), no. 2, 191--211.

\bibitem [MSI]{MSI}  \textit{Mirror Symmetry I}, Ed. S.-T. Yau, American Mathematical Society--International Press (1998).  This is the second edition of [EMM] above.

\bibitem [Per]{Per} U. Persson, \textit{Configurations of Kodaira Fibers on Rational Elliptic Surfaces}, Math. Z. 205 (1990), no. 1, 1--47. 

\bibitem [Pet]{Pet}  C. Peters, \textit{Monodromy and Picard-Fuchs Equations for Families of $K3$-Surfaces and Elliptic Curves}, Ann. Sci. \'{E}cole Norm. Sup. (4) 19 (1986), no. 4, 583--607. 

\bibitem  [PS]{PS} C. Peters and J. Stienstra, \textit{A Pencil of $K3$-Surfaces Related to Ap\'{e}ry's Recurrence for $\zeta(3)$ and Fermi Surfaces for Potential Zero}, Arithmetic of Complex Manifolds (Erlangen, 1988), 110--127, Lecture Notes in Math., 1399, Springer, Berlin (1989).


\bibitem [Sas]{Sas}  T. Sasai, \textit{Monodromy Representations of Homology of Certain Elliptic Surfaces}, J. Math. Soc. Japan, Vol. 26, No. 2, 296--305 (1974).

\bibitem [S-H]{S-H} Schmickler-Hirzebruch, \textit{Elliptische Fl\"{a}chen \"{u}ber $\mathbb{P}_1(\mathbb{C})$ mit drei Ausnahmefasern und die hypergeometrische Differentialgleichung}, Schriftenreihe des Mathematischen Instituts der Universit\"{a}t M\"{u}nster, Ser. 2, 33. Universit\"{a}t M\"{u}nster, Mathematisches Institut, M\"{u}nster (1985).

\bibitem [Shi]{Shi}  T. Shioda, \textit{On Elliptic Modular Surfaces}, J. Math. Soc. Japan, Vol.~24, 20--59 (1972). 


\bibitem [Sin]{Sin}  M. Singer, \textit{Algebraic Relations Among Solutions of Linear Differential Equations: Fano's Theorem}, Amer. J. Math., 110, no. 1, 115-143 (1988).

\bibitem [Sti]{Sti} P. Stiller, \textit{On the Uniformization of Certain Curves}, Pacific J. Math., 107, no. 1, 229--244 (1983). 

\bibitem [Var]{Var} V. Varadarajan, \textit{Meromorphic Differential Equations}, Exposition. Math. 9 (1991), no. 2, 97--188.

\bibitem [VY]{VY}  H. Verrill and N. Yui, \textit{Thompson Series, and the Mirror Maps of Pencils of $K3$ Surfaces}, contribution to this proceedings.


\bibitem [Voi]{Voi}  C. Voisin, \textit{Sym\'{e}trie Miroir}, Panoramas et Synth\`{e}ses, 2, Paris: Soci\'{e}t\'{e} Math\'{e}matique de France (1996).


\bibitem [Yos]{Yos} M. Yoshida, \textit{Fuchsian Differential Equations. With special emphasis on the Gauss-Schwarz theory}, Aspects of Mathematics, E11. Friedr. Vieweg \& Sohn, Braunschweig (1987).

\end{thebibliography}

\end{document}